\documentstyle[12pt]{article}

  \newcommand{\const}{\rm const}
  
  \newcommand{\supp}{\rm supp}

\begin{document}

   \begin{center}

{\bf A note about extension of functions} \\

\vspace{4mm}

{\bf  belonging to Sobolev - Grand Lebesgue Spaces.} \\

\vspace{6mm}

{\bf M.R.Formica, \ E.Ostrovsky, \ L.Sirota. }

 \end{center}

 \vspace{5mm}

\ Universit\`{a} degli Studi di Napoli Parthenope, via Generale Parisi 13, Palazzo Pacanowsky, 80132,
Napoli, Italy. \\
e-mail: \ mara.formica@uniparthenope.it \\

Department of Mathematics and Statistics, Bar-Ilan University,\\
59200, Ramat Gan, Israel. \\
e-mail: \ eugostrovsky@list.ru\\

Department of Mathematics and Statistics, Bar-Ilan University,\\
59200, Ramat Gan, Israel. \\
e-mail: \ sirota3@bezeqint.net \\

\vspace{5mm}

  \begin{center}

   {\bf Abstract}

  \end{center}

\vspace{4mm}

 \ We deduce an extension theorem for the so - called Sobolev - Grand Lebesgue Spaces defined on the  suitable subsets of
the  whole  finite - dimensional Euclidean space, and estimate the norms of correspondent extension operator, which
may be choosed as linear.\par

\vspace{5mm}

 \hspace{3mm} {\it  Key words and phrases.} Lebesgue - Riesz, Yudovich, ordinary Sobolev and Sobolev - Grand
 Lebesgue norm and spaces, Euclidean space, Lipschitz domain, ordinary and linear extension, linear bounded operator,  norm,
 measurable functions, semi - space, generating function, estimation.\par

\vspace{5mm}

\section{Introduction.  Previous results.}

\vspace{5mm}

 \hspace{3mm} Let $ \ B(G), \ $  where $ \ G \ $  is non - trivial subset of an Euclidean space
$  \  R^d: \ G \subset R^d \ $ be a family of Banach spaces defined on the class of functions defined  in turn on the
 support $ \ G; \ $ the norm in this space will be denoted $ \ ||f||B(G). \ $ \par
 \ It will be presumed henceforth that all the considered domains $ \ G \ $ are closures of  the non - empty open sets
and are Lipschitzian. The case when $ \ G = R^d \ $ is trivial for us and may be excluded. \par
 \ Put also for definiteness $ \  B = B(R^d). \ $ \par
 \  For instance, the classical Lebesgue - Riesz  $ \ L_p(G)  \ $ spaces equipped with the ordinary norm
$$
||f||_p(G) = ||f||L_p(G) := \left[ \ \int_G |f(x)|^p \ dx \ \right]^{1/p}, \ p \in [1,\infty), \ f: G \to R; \ x \in  G,
$$
as well as the famous Sobolev spaces $ \ W^m_p(G), \ m = 1,2,\ldots: \ $
\begin{equation} \label{Sobolev space}
||f||W^m_p(G) \stackrel{def}{=}  \max_{\alpha, \ |\alpha| \le m} ||D^{\alpha} f ||L_p(G),
\end{equation}

$  \ ||f||W^m_p := ||f||W^m_p(R^d); \ $ see e.g.  \cite{Adams}, \cite{Leoni}, \cite{Mazja 1}, \cite{Mazja 2}, \cite{Sobolev}, \cite{Taylor}. \par

\ Here as ordinary

$$
  \alpha  = \vec{\alpha} = (\alpha_1, \alpha_2, \ldots, \alpha_d); \hspace{3mm}  \alpha_j = 0,1,2,\ldots;
$$

$$
|\alpha| := \sum_{j=1}^d \alpha_j; \hspace{3mm} D^{\alpha} f := \frac{D^{|\alpha|} f}{\partial x_1^{\alpha_1} \partial x_2^{\alpha_2} \ldots \partial x_d^{\alpha_d} }
$$
and all the derivatives are understood in the weak (Sobolev) sense.  Of course,  $ \ D^0 f = f. \ $ \par

\vspace{3mm}

 \ Let the function $ \ f: G \to R \ $ belongs to some space $ \ B(G). \ $ By definition, the function $ \  \tilde{f}, \ R^d \to R \ $
is named as an {\it extension}  of the function $ \ f \ $ (from the set $ \ G),  \ $  iff

$$
\forall x \in G \ \Rightarrow \tilde{f}(x) = f(x)
$$
and wherein  $ \ ||\tilde{f}|| B < \infty.  \ $  \par
 \ If there exists a {\it linear  bounded} operator $ \  L: \ \tilde{f} = L f = L_G f, \ f \in B(G),  \ $ where again  $ \ \tilde{f} \ $ is an extension for $ \ f, \ $
 such that

$$
K = K(G) = K(B,G) \stackrel{def}{=} \sup_{0 \ne f \in B(G)} \left\{ \  \frac{||\ Lf \ ||B \ }{ ||f||B(G)} \ \right\} < \infty,
$$
then the extension $ \ \tilde{f} = Lf  = L_G f \ $ is named linear. \par

\vspace{3mm}

 \ It is known, see e.g. \cite{Burenkov 1}, \cite{Burenkov 2}, \cite{Fefferman}, \cite{Taylor} that under imposed conditions on the Lipschitz
domain $ \ G \ $ and for  the Sobolev's spaces $ \ B(G) = W^m_p(G), \ m \ge 1, \ p \ge 1 \ $ (the case $ \ m = 0 \ $ is trivial) \
 there exist a linear extension operator. See also the works \cite{Garcia}, \cite{Koskela},  \cite{Mazja 1}, \cite{Mazja 2}, \cite{Romanov}  etc.\par

\vspace{4mm}

 \hspace{3mm} {\bf We will ground in offered report  that this linear  operator there exists also for the so - called
 Sobolev - Grand Lebesgue Spaces. } \par

\vspace{5mm}

\section{Main result.}

\vspace{5mm}

\begin{center}

 \  {\sc Sobolev -  Grand Lebesgue Spaces.} \par

\end{center}

\vspace{4mm}

 \hspace{3mm} Let $ \  (a,b) = \const, 1 \le a < b \le \infty.  \  $ Let also $ \ \psi = \psi(p), \ p \in (a,b) \ $  be certain numerical valued
 measurable strictly positive: $ \  \inf_{p \in (a,b)} \psi(p) > 0   \ $  function, not necessary to be bounded. Denotation:

 $$
  (a,b) := \supp (\psi); \  \Psi[a,b] := \{ \  \psi: \ \supp (\psi) = (a,b)  \ \},
 $$

$$
\Psi \stackrel{def}{=} \cup_{ (a,b): \ 1 \le a < b \le \infty } \Psi[a,b].
$$

\vspace{4mm}

\hspace{3mm} {\bf Definition 2.1.}, see e.g. \cite{Ostrovsky Sob GLS}. \
 \ The Sobolev - Grand Lebesgue Space $ \ S[m, G, \psi] \ $
 based on the set $ \ G, \ G \subset R^d \ $ is defined as a set of all  (measurable) functions  having a finite norm

\vspace{3mm}

\begin{equation} \label{def SGLS norm}
||f||S[m, G, \psi]  \stackrel{def}{=}  \sup_{p \in (a,b)} \left\{ \ \frac{||f||W^m_p(G)}{\psi(p)} \ \right\}.
\end{equation}

\vspace{4mm}

 \ The particular case $ \ m = 0, \ $ i.e. when

\begin{equation} \label{Gpsi}
||f||G\psi = ||f||G\psi(a,b)  \stackrel{def}{=} ||f|||S[0, G, \psi]  = \sup_{p \in (a,b)} \left\{ \ \frac{||f||L_p(G)}{\psi(p)} \ \right\}
\end{equation}
and under  some additional restrictions on the  {\it generating function} $ \ \psi = \psi(p) \ $
correspondent to the so - called {\it Yudovich spaces,}  see  \cite{Yudovich 1}, \cite{Yudovich 2}.
These spaces was applied at first in the theory of Partial Differential Equations (PDE),  see  \cite{Chen},  \cite{Crippa}. \par

\vspace{3mm}

 \ A general case of these spaces when $ \ m = 0 \ $ are named as
 the  classical Grand Lebesgue Spaces (GLS)  $ \ G \psi, \ \psi \in \Psi. \ $ These spaces are
 investigated in many works, see e.g. \cite{ErOs}, \cite{Ermakov}, \cite{Fiorenza2}, \cite{Fiorenza-Formica-Gogatishvili-DEA2018},
\cite{fioforgogakoparakoNAtoappear}, \cite{fioformicarakodie2017}, \cite{formicagiovamjom2015}, \cite{Formica Ostrovsky Sirota weak dep},
\cite{KosOs}, \cite{KozOsSir2017}, \cite{KosOs equivalence}, \cite{Liflyand}, \cite{Ostrovsky1}.  The general case of Sobolev -
Grand Lebesgue Spaces appears at first perhaps in the article \cite{Ostrovsky Sob GLS}, where was investigated the modulus of
continuity of the functions belonging to these spaces. \par

\vspace{4mm}

\ {\bf Theorem 2.1.} Assume that all the  formulated before restrictions are satisfied,  indeed: that the Lipschitz  domain $ \ G \ $ is closure
of the non - empty open subset
of whole Euclidean space $ \ R^d. \ $  We propose that for arbitrary Sobolev - Grand Lebesgue Space $ \ S[m, G, \psi] \ $
there  exists a {\it linear  bounded} extension  operator  $ \ L = L_G. \ $  \par

\vspace{4mm}

 \hspace{3mm}  {\bf Proof.} One can suppose $ \ d \ge 2 \ $ and that $ \ G =  \vec{x} = \{x_j\} = \{x_1, x_2, \ldots, x_{d-1}, x_d \}, \ $  where
$ \ \vec{x} = (\tilde{x}, x_d); \  \tilde{x} = \{ \ x_1, x_2, \ldots, x_{d-1} \ \}; \ $ so that
  $ \ \vec{x} \in G \ \Leftrightarrow  x_d \ge 0; \ $ on the other words, upper semi - space;
  see e.g. \cite{Burenkov 1}, \cite{Burenkov 2}, \cite{Fefferman}. Let us define the following
 extension operator  $ \ Lf(x) := f(x), \  x \in G, \  f(\cdot) \in S[m,G,\psi]; \ $

 $$
 Lf(x): = \sum_{k=1}^{d + 1} c_k f(\tilde{x}, \ - k x_d), \ x_d < 0.
 $$
 \ The coefficients $ \ \{ c_k\} \ $ may be uniquely determined from the following system of linear equations

 $$
 \sum_{k=1}^{m+1} (-k)^l \ c_k = 1; \ l =  0,1,\ldots,d.
 $$
  \ Suppose that $ \ f \in S[m,G,\psi]; \ $ one can assume without loss of generality $ \ ||f||S[m,G,\psi] = 1. \ $  Then for all the values $ \ p \in (a,b) \ $

 $$
||f||W^m_p \le \psi(p), \ p \in (a,b),
 $$
therefore

$$
\forall p\in (a,b), \hspace{3mm}  \forall \alpha: |\alpha| \le m \ \Rightarrow ||D^{\alpha} f||_p(G)  \le \psi(p).
$$
 \ Introduce the functions

$$
g_k(\tilde{x},y) := f( \tilde{x}, - k y), \ y \le 0;
$$
then

$$
D^{\alpha} g_k = (-k)^{\alpha_d} \ f^{(\alpha)}( \tilde{x}, - k y),
$$

$$
||D^{\alpha} g_k||_p(G) = k^{\alpha_d - 1/p} ||D^{\alpha}  f||_p(G)  \le  k^{\alpha_d} ||D^{\alpha} f||_p(G),
$$
therefore

$$
\forall \alpha: |\alpha| \le m \ \Rightarrow \sum_{k=1}^{m+1} |c_k| \ k^m \cdot ||D^{\alpha} f||_p(G) \le
$$

$$
\sum_{k=1}^{m+1} |c_k| \ k^m \cdot \psi(p), \ p \in (a,b).
$$

 \  Following,  by virtue of triangle inequality for Lebesgue - Riesz spaces

$$
||L[f]||W^m_p(R^d \setminus G) \le C(d,m) \ \psi(p), \ C(d,m) < \infty,
$$

$$
||Lf||S[m,R^d,\psi]  \le 1 + C(d,m)  < \infty,
$$
Q.E.D. \par

\vspace{5mm}

\section{Concluding remarks.}

\vspace{5mm}

 \hspace{3mm} It is interest in our opinion to compute the exact value of extension constant for Sobolev - Grand Lebesgue Spaces,
 as well as to generalize the extension theorem on the anisotropic spaces. \par

\vspace{6mm}

\vspace{0.5cm} \emph{Acknowledgement.} {\footnotesize The first
author has been partially supported by the Gruppo Nazionale per
l'Analisi Matematica, la Probabilit\`a e le loro Applicazioni
(GNAMPA) of the Istituto Nazionale di Alta Matematica (INdAM) and by
Universit\`a degli Studi di Napoli Parthenope through the project
\lq\lq sostegno alla Ricerca individuale\rq\rq .\par


\vspace{5mm}




 \begin{thebibliography}{44}

\bibitem{Adams}
{\bf  Adams, R.A.} (1975), {\it Sobolev spaces.}  Academic Press, New York.

\bibitem{Buldygin}
{\bf V.V. Buldygin V.V., D.I.Mushtary, E.I.Ostrovsky, M.I.Pushalsky.} {\it New Trends in Probability Theory and Statistics.}
Mokslas, (1992), V.1, p. 78 \ - \ 92; Amsterdam, Utrecht, New York, Tokyo.

\bibitem{Burenkov 1}
{\bf  Burenkov V. I.} {\it On a certain method for extending differentiable functions.} (Russian). Trudy Mat. Inst. Steklov 140 (1976), 27 \ - \ 67;
 English translation.: Proc. of Steklov Inst. Math.; {\bf 140,}  (1976).

\bibitem{Burenkov 2}
{\bf  Burenkov V. I.} {\it On the extension of functions with preservation of semi-norm.} (Russian). Dokl. Akad. Nauk SSSR, {\bf 228,}  (1976), 779 \ - \ 782;
 English translation.: Soviet Math. Dokl. {\bf 17, }  (1976). MR 0412796

\bibitem{Capone1}
{\bf Capone C, Formica M.R, Giova R.} {\it Grand Lebesgue spaces with respect to measurable functions.}
Nonlinear Analysis 2013; 85: 125 \ - \ 131.

\bibitem{Capone2}
{\bf Capone C, and Fiorenza A.} {\it On small Lebesgue spaces. Journal of function spaces and applications.}
2005; 3; \ 73 \ - \ 89.

\bibitem{Chen}
{\bf Qionglei Chen, Changxing Miao and Xiaoxin  Zheng.}  {\it The two - dimensional Euler equation in Yudovich and bmo - type spaces.} \\
arXiv:1311.0934v4 [math.AP] 10 Jan 2019


\bibitem{Crippa}
{\bf Gianluca Crippa and Giorgio Stefani.} {\it  An elementary proof of existence and uniqueness for the Euler flow in localized
Yudovich spaces.} \\
arXiv:2110.15648v1 [math.AP] 29 Oct 2021


\bibitem{ErOs}
{\bf S. V. Ermakov, E. I. Ostrovsky.}  {\it Central limit theorem for weakly dependent Banach \ - \ space valued random variables.} Theory Probab. Appl.,
{\bf 30}, 2, (1986), 391 \ - \ 394.

\bibitem{Ermakov}
{\bf S.V.Ermakov, and E. I. Ostrovsky.} {\it Continuity Conditions,  Exponential Estimates, and the Central Limit Theorem for Random Fields.}
Moscow, VINITY, (1986), (in Russian).

\bibitem{Fefferman}
{\bf Charles L.Fefferman, Arie Israel, and Garving K. Luli.}   {\it Sobolev extension by linear operators.}
Dorderecht,  (2012).




\bibitem{Fiorenza2}
{\bf Fiorenza A., and Karadzhov G.E.} {\it Grand and small Lebesgue spaces and their analogs.} Consiglio Nationale Delle Ricerche,
Instituto per le Applicazioni del Calcoto Mauro Picone,
Sezione di Napoli, Rapporto tecnicon. 272/03, (2005).

\bibitem{Fiorenza-Formica-Gogatishvili-DEA2018}
{\bf A.~Fiorenza, M.~R.~Formica} and {\bf A.~Gogatishvili.} {\it On
grand and small Lebesgue and Sobolev spaces and some applications to
PDE's}. \emph{Differ. Equ. Appl.} \textbf{10} (2018), no.~1, 21--46.

\bibitem{fioforgogakoparakoNAtoappear}
{\bf A.~Fiorenza, M. R.~Formica, A.~Gogatishvili, T.~Kopaliani} and
{\bf J.~M. Rakotoson.} {\it Characterization of interpolation
between grand, small or classical Lebesgue spaces}. Preprint
arXiv:1709.05892, Nonlinear Anal., {to appear}.

\bibitem{fioformicarakodie2017}
{\bf A.~Fiorenza, M.~R.~Formica} and {\bf J.~M. Rakotoson.} {\it
Pointwise estimates for {$G\Gamma$}-functions and applications}.
Differential Integral Equations \textbf{30} (2017), no.~11-12,
809--824.

\bibitem{formicagiovamjom2015}
{\bf M.~R. Formica} and {\bf R.~Giova.} {\it Boyd indices in
generalized grand Lebesgue spaces and applications}. Mediterr. J.
Math. \textbf{12} (2015), no.~3, 987--995.

\bibitem{Formica Ostrovsky Sirota weak dep}
{\bf M.R.Formica, E.Ostrovsky, L.Sirota} {\it Central Limit Theorem in Lebesgue-Riesz spaces
for weakly dependent random sequences.} \\
arXiv:1912.00338v2 [math.PR] 3 Dec 2019


\bibitem{Garcia}
{\bf Miguel Garcia, Ia - Bravo, Tapio Rajala, and Jurki Takanen.} {\it Two - sides boundary points of Sobolev - extension
domains on Euclidean spaces.} \\
arXiv:2111.01079v1  [math.CA]  1 Nov 2021


\bibitem{Koskela}
{\bf P. Koskela, D.Yang, and Y. Zhou.}  {\it A Jordan Sobolev extension domain.} Ann. Acad. Sci. Fenn. Math.  {\bf 35,}
(2010), no. 1, 309 \ - \ 320.




\bibitem{KosOs}
{\bf Yu.V. Kozachenko and E.I. Ostrovsky.} {\it Banach spaces of random variables of subgaussian type.} Theory Probab. Math. Stat., Kiev, (1985), v. 43,
\ 42 \ - \ 56, (in Russian).

\bibitem{KozOsSir2017}
{\bf Kozachenko Yu.V., Ostrovsky E., Sirota L.}
{\it Relations between exponential tails, moments and moment generating functions for random variables and vectors.} \\
arXiv:1701.01901v1 [math.FA] 8 Jan 2017

\bibitem{KosOs equivalence}
{\bf Kozachenko Yu.V., Ostrovsky E., Sirota L.}  {\it Equivalence between tails, Grand Lebesgue Spaces and Orlicz norms for random variables without
Kramer's condition.} Bulletin of KSU, Kiev, 2018, {\bf  4,} pp. 20 - 29.

\bibitem{Liflyand}
{\bf E.Liflyand,  E.Ostrovsky,  L.Sirota.} {\it Structural Properties of Bilateral Grand Lebesgue Spaces.} Turk. J. Math.; {\bf 34}, (2010), 207 \ - \ 219.


\bibitem{Leoni}
{\bf Leoni, Giovanni.}  (2009).   {\it A First Course in Sobolev Spaces.} Graduate Studies in Mathematics, {\bf 105.}
 American Mathematical Society. pp. xvi + 607.\\
 ISBN 978-0-8218-4768-8. MR 2527916. Zbl 1180.46001.





\bibitem{Mazja 1}
{ \bf V. Maz'ya.} {\it Extension of functions from Sobolev spaces.}  Zap. Nauchn. Sem. Leningrad. Otdel. Mat. Inst.
Steklov. (LOMI), {\bf  113}, (1981), 231 \ - \ 236.


\bibitem{Mazja 2}
{\bf Maz'ja, Vladimir G.}  (1985).  {\it Sobolev Spaces.}  Springer Series in Soviet Mathematics, Berlin–Heidelberg–New York: Springer-Verlag, pp. xix+486,\\
 doi:10.1007/978-3-662-09922-3, ISBN 0-387-13589-8, MR 0817985, Zbl 0692.46023



\bibitem{Ostrovsky1}
{\bf E.I. Ostrovsky.} {\it Exponential Estimations for Random Fields.} Moscow \ - \ Obninsk, OINPE, (1999), in Russian.


\bibitem{Ostrovsky Sob GLS}
{\bf Ostrovsky E., Sirota L.}   {\it Module of continuity for the functions belonging to the Sobolev-Grand Lebesgue Spaces.}\\
arXiv:1006.4177v1  [math.FA]  21 Jun 2010


\bibitem{Romanov}
{\bf A.S.Romanov.}  {\it On the extension of functions that belong to Sobolev spaces.} Sibirsk. Mat. Zh., {\bf 34},  (1993),
149 \ - \ 152; English transl. in Siberian Math. J. {\bf 34},  (1993), 723 \ - \ 726.


\bibitem{Sobolev}
{\bf Sobolev, S.L.} (1963). {\it Some applications of functional analysis in mathematical physics.}  Amer. Math. Soc.



\bibitem{Taylor}
{\bf Steve Taylor.} {\it An introduction to Sobolev spaces.}  (2015), American Mathematical Society. pp. xvi + 769.




\bibitem{Yudovich 1}
{\bf V. I. Yudovich.} {\it Nonstationary flow of an ideal incompressible liquid.}  Zh. Vych. Mat., {\bf 3,}  (1963), 1032  \ - \ 1066.




\bibitem{Yudovich 2}
{\bf  V. I. Yudovich.} {\it Uniqueness theorem for the basic nonstationary problem in the dynamics of an ideal
incompressible fluid.}  Mathematical Research Letters, {\bf 2.} (1995), 27 \ - \ 38.


\end{thebibliography}
\end{document}